\documentclass[12pt, final, reqno]{amsart}
\usepackage{amssymb}

\theoremstyle{plain}
\newtheorem{theorem}{Theorem}
\newtheorem{lemma}{Lemma}
\newtheorem{corollary}[theorem]{Corollary}

\newcommand\BP{B(H)^+}
\newcommand\BS{B_s(H)}
\newcommand\Mnp{M_n(\mathbb C)^+}
\newcommand\N{\mathbb N}
\newcommand\R{\mathbb R}
\newcommand\rank{\operatorname{rank}}
\newcommand\rng{\operatorname{rng}}

\begin{document}
\vglue -20pt

\centerline{\large{\textbf{ORDER-AUTOMORPHISMS OF THE SET OF}}}

\centerline{\large{\textbf{BOUNDED OBSERVABLES}}}

\vskip 15pt

\centerline{LAJOS MOLN\'AR}

\vskip -5pt

\centerline{Institute of Mathematics and Informatics}

\vskip -5pt

\centerline{University of Debrecen}

\vskip -5pt

\centerline{H-4010 Debrecen, P.O.Box 12}

\vskip -5pt

\centerline{Hungary}

\vskip 10pt

\centerline{phones: dept. ++ 52 512 900/2815;
                    fax ++ 52 416 857}

\vskip 5pt

\centerline{e-mail: \texttt{molnarl@math.klte.hu}, }

\vskip 10pt

\centerline{Running title:}

\centerline{\scshape{AUTOMORPHISMS OF THE SET OF BOUNDED
OBSERVABLES}}

\pagestyle{myheadings} \markboth{\textsc{\SMALL L. Moln\'ar}}
{\textsc{\SMALL Automorphisms of the set of bounded observables}}

\normalsize \vskip 20pt

\centerline{\textsc{Abstract}} Let $H$ be a complex Hilbert space
and denote by $\BS$ the set of all self-adjoint bounded linear
operators on $H$. In this paper we describe the form of all
bijective maps (no linearity or continuity is assumed) on $\BS$
which preserve the order $\leq$ in both directions.

\newpage
\section{Introduction and statements of the results}

In the Hilbert space framework of quantum mechanics the bounded
observables are represented by self-adjoint bounded linear
operators. If $H$ denotes the underlying Hilbert space, then these
operators form the set $\BS$ on which usually several operations
and relations are considered. The automorphisms of $\BS$ with
respect to those operations and/or relations are, just as with any
algebraic structure in mathematics, of remarkable importance.

First of all, $\BS$ with the usual addition, scalar multiplication
and Jordan product forms a Jordan algebra. It is a well-known
result that the corresponding automorphisms of $\BS$ are
implemented by unitary or antiunitary operators of $H$ (see, for
example, Ref. \cite{CVLL}, where the automorphisms of some other
important structures appearing in the probabilistic aspects of
quantum mechanics are also treated).

The aim of this paper is to determine another class of
automorphisms of $\BS$. Namely, we equip the set $\BS$ with the
usual order among self-adjoint operators. That is, for any $A,B\in
\BS$, we write $A\leq B$ if $\langle Ax,x\rangle\leq \langle B
x,x\rangle$ holds for every $x\in H$. Alternatively, in the
language of quantum mechanics, the bounded observable $A$ is said
to be less than or equal to the bounded observable $B$ if the
expected value of $A$ in any state is less than or equal to the
expected value of $B$ in the same state. The relation $\leq$ is no
doubt an important one among observables.

In what follows we determine all the automorphisms of $\BS$ as a
partially ordered set with the relation $\leq$ (this is done in
the main result of the paper Theorem~\ref{T:lj2}) and also present
some corollaries (Corollary~\ref{C:lj3}, Corollary~\ref{C:lj4})
that we believe have physical meaning.

We begin with the following proposition on which the proof of our
main result rests. Let $H$ be a complex Hilbert space and let
$\BP$ denote the cone of all positive operators on $H$ (that is, the set
of all $A\in \BS$ for which $\langle Ax,x\rangle \geq 0$ holds
for every $x\in H$). Our first result describes the form of all
bijective
maps on $\BP$ which preserve the order $\leq$ in both directions.

\begin{theorem}\label{T:lj1}
Assume that $\dim H>1$. Let $\phi: \BP\to \BP$ be a bijective map
with the property that
\[
A\leq B \Longleftrightarrow \phi(A)\leq \phi(B)
\]
holds whenever $A,B\in \BP$. Then there exists an invertible bounded
either linear or conjugate-linear operator $T:H \to H$ such that
$\phi$ is of the form
\[
\phi(A)=TAT^* \qquad (A\in \BP).
\]
\end{theorem}

After having proved this result, it will be easy to deduce the main
result of the paper that follows.

\begin{theorem}\label{T:lj2}
Suppose that $\dim H>1$. Let $\phi: \BS\to \BS$ be a bijective map
with the property that
\[
A\leq B \Longleftrightarrow \phi(A)\leq \phi(B)
\]
holds whenever $A,B\in \BS$. Then there exists an operator $X\in
\BS$ and an invertible bounded either linear or conjugate-linear
operator $T:H \to H$ such that $\phi$ is of the form
\[
\phi(A)=TAT^* +X\qquad (A\in \BS).
\]
\end{theorem}

This result has some corollaries that seem worth mentioning. In
the first one we determine the form of all bijective
transformations on $\BS$ which preserve the order and the
commutativity in both directions (in quantum mechanics, instead of
commutativity they usually use the word compatibility for this
important concept).

\begin{corollary}\label{C:lj3}
Assume that $\dim H>1$. Let $\phi:\BS\to \BS$ be a bijective map
which preserves the order and the commutativity in both
directions. Then there is an either unitary or antiunitary
operator $U:H\to H$, a positive scalar $\lambda$, and a real
number $\mu$ such that $\phi$ is of the form
\[
\phi(A)=\lambda UAU^* +\mu I \qquad (A\in \BS).
\]
\end{corollary}

The next corollary describes all the bijective maps on $\BS$ which
preserve the order and the complementarity in both directions (two
observables are called complementary if the range of any
nontrivial projection from the range of the spectral measure of
the first observable has zero intersection with the range of any
nontrivial projection from the range of the spectral measure of
the second observable). Although this latter concept is in some
sense opposite to compatibility, as it turns out below we still
have the same form for $\phi$ as above.

\begin{corollary}\label{C:lj4}
Suppose that $\dim H>1$. Let $\phi:\BS\to \BS$ be a bijective map
which preserves the order and the complementarity in both
directions. Then there is an either unitary or antiunitary
operator $U:H\to H$, a positive scalar $\lambda$, and a real
number $\mu$ such that $\phi$ is of the form
\[
\phi(A)=\lambda UAU^* +\mu I \qquad (A\in \BS).
\]
\end{corollary}

Finally, our last corollary characterizes those bijective maps on
$\BS$ which preserve the order and the orthogonality in both
directions (two operators $A,B \in \BS$ are called orthogonal if
$AB=0$ which is just equivalent to the mutual orthogonality of the
ranges of $A$ and $B$).

\begin{corollary}\label{C:lj5}
Assume that $\dim H>1$. Let $\phi:\BS\to \BS$ be a bijective map
which preserves the order and the orthogonality in both
directions. Then there is an either unitary or antiunitary
operator $U:H\to H$, and a positive scalar $\lambda$ such that
$\phi$ is of the form
\[
\phi(A)=\lambda UAU^* \qquad (A\in \BS).
\]
\end{corollary}

Closing this section we note that all the above statements could
be conversed, that is, if a map $\phi$ is of any of the above
forms, then it necessarily preserves the corresponding properties.

\section{Proofs}

This section is devoted to the proofs of our results. We begin
with the following auxiliary results. If $A$ is a bounded linear
operator, then $\rng A$ denotes its range. The rank of $A$ is, by
definition, the algebraic dimension of $\rng A$ and it is denoted
by $\rank A$. For any $A\in \BP$, $\sqrt A$ stands for the unique
positive linear operator whose square is $A$.

\begin{lemma}\label{L:rng}
Let $A,B \in \BP$ be such that $\rank A=1, \rank B< \infty$. We
have $\lambda A\leq B$ for some positive scalar $\lambda$ if and
only if ${\rng A}\subset {\rng B}$.
\end{lemma}

\begin{proof}
We recall the following nice result of Busch and Gudder (Ref.
\cite[Theorem 3]{BuschGudder}): if $B\in \BP$, $x\in H$, and $P$
is the rank-1 projection projecting onto the subspace generated by
$x$, then we have $\lambda P\leq B$ for some positive scalar
$\lambda$ if and only if $x$ is in the range of $\sqrt B$. As in
our case $B$ is a finite rank operator, it follows from the
spectral theorem that $\rng B=\rng \sqrt B$. Since the positive
rank-1 operators are exactly the positive scalar multiples of
rank-1 projections, we obtain the assertion.
\end{proof}

\begin{lemma}\label{L:>n+1}
Let $A\in \BP$ and $n\in \N$. We have $\rank A>n+1$ if and only if there
are operators $E,F\in \BP$ such that $E,F\leq A$, $\rank E=n, \rank F>1$
and there is no $G\in \BP$ of rank 1 with $G\leq E,F$.
\end{lemma}

\begin{proof}
Suppose that $\rank A>n+1$. We assert that there exists a finite
rank operator $A'\in \BP$ such that $A'\leq A$ and $\rank A'>n+1$.
In case $A$ is of finite rank, this is trivial. If $A$ is compact
and not of finite rank, then by the spectral theorem of compact
self-adjoint operators we can verify our claim very easily.
Finally, if $A$ is non-compact, then using the spectral theorem of
self-adjoint operators and the properties of the spectral
integral, we can find an infinite rank projection $P$ on $H$ and a
positive scalar $\lambda$ such that $\lambda P\leq A$ from which
the existence of an appropriate operator $A'$ follows.

Clearly, $A'$ can be written as the sum of positive scalar multiples
of pairwise orthogonal rank-1 projections. Let $E$ be the sum of the
first $n$ terms in this sum and let $F$ be the sum of the remaining
part.
It is easy to see that $E,F$ have the required property.
In fact, the non-existence of $G$ follows from Lemma~\ref{L:rng}.

To prove the converse, suppose that there are operators $E,F\in
\BP$ with the properties formulated in the lemma. It follows from
the relation $E,F\leq A$ that $E,F$ are of finite rank and $\rng
E, \rng F \subset \rng A$. As there is no positive rank-1 operator
$G$ with $G\leq E,F$, by Lemma~\ref{L:rng} we have $\rng E\cap
\rng F=\{ 0\}$. So, $\rng A$ contains two subspaces with trivial
intersection the sum of whose dimensions is greater than $n+1$.
This shows that $\rank A> n+1$, completing the proof of the lemma.
\end{proof}

Now, we are in a position to prove our first theorem.

\begin{proof}[Proof of Theorem~\ref{T:lj1}]
We first remark that our proof is based on a beautiful result of
Rothaus (Ref. \cite{Rothaus}) concerning the automatic linearity
of bijective maps between closed convex cones in normed spaces
preserving order in both directions. In the paper Ref.
\cite{Rothaus} conclusions of that kind were reached under some
quite restrictive assumptions. In our present situation, that is,
when the normed space in question is an operator algebra, those
assumptions are fulfilled exactly when the underlying Hilbert
space $H$ is finite dimensional. Accordingly, the main point of
our proof is to reduce the problem to the finite dimensional case.
This is in fact what we are going to do below.

Clearly, $\phi(0)=0$. We prove that $\phi$ preserves the rank of
operators. In fact, we show that the assertion that
\[
\rank A=k \Longleftrightarrow \rank \phi(A)=k
\]
$(k=1,\ldots, n)$ holds for every $n\in \N$. To begin, as for the
case $n=1$, we remark that a nonzero operator $A\in \BP$ is of
rank 1 if and only if the operator interval $[0,A]$ is total under
the partial ordering $\leq$, that is, every two elements of it are
comparable. Suppose that our assertion is true for some $n\in \N$.
We show that in that case it holds also for $n+1$. Let $A\in \BP$
be of rank $n+1$. By our assumption of induction, it follows that
the rank of $\phi(A)$ is at least $n+1$. Suppose that $\rank
\phi(A)> n+1$. Using Lemma~\ref{L:>n+1} and the order preserving
property of $\phi$ we obtain that $\rank A>n+1$ which is a
contradiction. Therefore, we have $\rank \phi(A)=n+1$. Referring
to the fact that $\phi^{-1}$ shares the same properties as $\phi$,
we obtain the desired assertion.

We now prove that if $A_1, \ldots, A_n\in \BP$ are of rank 1, then
their ranges are linearly independent if and only if so are the
ranges of $\phi(A_1), \ldots, \phi(A_n)$. (A system of
1-dimensional subspaces in $H$ of $n$ members is called linearly
independent if they cannot be included in an $(n-1)$-dimensional
subspace.) This statement is clear for $n=1$. Suppose that it
holds for $n$ and prove that it then necessarily holds also for
$n+1$. Let $A_1, \ldots ,A_n,A_{n+1}$ be rank-1 operators with
linearly independent ranges and assume that this is not the case
with the ranges of $\phi(A_1), \ldots, \phi(A_n), \phi(A_{n+1})$.
Then these ranges can be included in an at most $n$-dimensional
subspace implying that there is a rank-$n$ operator $B\in \BP$
such that $\phi(A_1), \ldots, \phi(A_{n+1})\leq B$. By the
rank-preserving property of $\phi$ we have a rank-$n$ operator
$A\in \BP$ such that $A_1, \ldots A_n, A_{n+1} \leq A$. By
Lemma~\ref{L:rng} this implies that $\rng A_1, \ldots, \rng
A_{n+1} \subset \rng A$ and it follows that the ranges of $A_1,
\ldots ,A_{n+1}$ can be included in an $n$-dimensional subspace of
$H$ which is a contradiction. This verifies our claim.

Fix rank-1 operators $A_1, \ldots, A_n\in \BP$ with linearly
independent ranges which generate the $n$-dimensional subspace
$H_n$ of $H$. Denote by $H_n'$ the $n$-dimensional subspace of $H$
generated by the ranges of $\phi(A_1), \ldots, \phi(A_n)$. We
assert that an operator $T\in \BP$ acts on $H_n$ if and only if
$\phi(T)$ acts on $H_n'$. (We say that an operator $T$ acts on the
closed subspace $M$ of $H$ if $M$ is an invariant subspace of $T$
and $T$ is zero on the orthogonal complement of $M$.) This will
follow from the following observation: the positive finite rank
operator $T$ acts on $H_n$ if and only if for every rank-1
operator $A$ for which the ranges of $A_1, \ldots ,A_n, A$ are
linearly independent we have $A\nleq T$. To see this, suppose that
$T$ acts on $H_n$. If $A\leq T$, then we have $\rng A\subset \rng
T\subset H_n$ implying that the ranges of $A_1, \ldots ,A_n, A$
cannot be linearly independent. This gives us the necessity. As
for the sufficiency, suppose that $T$ does not act on $H_n$. Then
there exists a unit vector $x$ in the range of $T$ which does not
belong to $H_n$. On the other hand, as $x\in \rng T$, by
Lemma~\ref{L:rng} it follows that a positive scalar multiple of
the rank-1 projection projecting onto the subspace generated by
$x$ is less than or equal to $T$. This gives us a rank-1 operator
$A$ for which the ranges of $A_1, \ldots, A_n, A$ are linearly
independent and we have $A\leq T$. This proves our claim.

So, for any $n$-dimensional subspace $H_n$ of $H$, there exists an
$n$-dimensional subspace $H_n'$ of $H$ such that for every $T\in
\BP$, $T$ acts on $H_n$ if and only if $\phi(T)$ acts on $H_n'$.
This gives rise to a bijective transformation $\psi$ on the cone
$\Mnp$ of all positive $n\times n$ complex matrices which
preserves the order in both direction. (Here positivity is used in
the operator theoretical sense, so our concept of positivity is
just the same as positive semidefiniteness in matrix theory.)

Since $\phi$ preserves the rank, it follows that $\psi$ preserves
the rank-$n$ matrices in both directions. The set of all such
matrices is just the interior of $\Mnp$ in the real normed space
of all $n\times n$ Hermitian matrices. Now, the result Ref.
\cite[Proposition 2]{Rothaus} of Rothaus on the linearity of order
preserving maps can be applied and it gives us that $\psi$ is
linear on the set of all rank-$n$ elements in $\Mnp$. We show that
$\psi$ is linear on the whole set $\Mnp$. Pick $A,B\in \Mnp$. Then
there are sequences $(A_k), (B_k)$ of rank-$n$ elements in $\Mnp$
which are monotone decreasing with respect to the order $\leq$ and
$A_k \to A$, $B_k\to B$. It is clear that the equalities $A=\inf_k
A_k$, $B=\inf_k B_k$ and $A+B=\inf_k (A_k+B_k)$ hold in the
partially ordered set $\Mnp$. By the order preserving property of
$\psi$ we obtain that $\psi(A)=\inf_k \psi(A_k)$, $\psi(B)=\inf_k
\psi(B_k)$ and $\psi(A+B)=\inf_k \psi(A_k+B_k)$. The sequences
$\psi(A_k)$, $\psi(B_k)$, $\psi(A_k+B_k)$ are monotone decreasing
and bounded below. By Vigier's theorem (Ref. \cite[4.1.1.
Theorem]{Murphy}) they necessarily converge (strongly) to their
infima. Now, by the partial additivity property of $\psi$ which
has been obtained above as a consequence of Rothaus's result, we
have
\[
\psi(A+B)=\lim_k \psi(A_k+B_k)=\lim_k \psi(A_k)+\lim_k
\psi(B_k)=\psi(A)+\psi(B).
\]
So, $\psi$ is additive on $\Mnp$ and one can prove in the same way
that it is positive homogeneous as well. Since every pair of
finite rank elements in $\BP$ can be embedded into a matrix space
$\Mnp$, we deduce that $\phi$ is additive and positive homogeneous
on the set of all finite rank elements in $\BP$.

Since every finite sum $\sum_i \lambda_i P_i$, where the
$\lambda_i$'s are positive numbers and the $P_i$'s are projections
of not necessarily finite rank, is the strong limit of a monotone
increasing net of finite rank elements in $\BP$, one can prove in
a very similar way as above that $\phi$ is additive and positive
homogeneous on the set of all such finite sums. Finally, using the
fact that every operator in $\BP$ is the norm limit of a monotone
increasing sequence of operators of the form $\sum_i \lambda_i
P_i$ (this follows form the spectral theorem), repeating the above
argument once again, we obtain that $\phi$ is additive and
positive homogeneous.

Extend $\phi$ from $\BP$ to $B(H)_s$ in the obvious way, that is,
define $\tilde \phi(T)=\phi(A)-\phi(B)$ for every $T\in \BS$ and
$A,B\in \BP$ for which $T=A-B$. It is easy to check that $\tilde
\phi:B(H)_s \to B(H)_s$ is a linear transformation which preserves
the order in both directions. To see the less trivial part of this
last assertion, suppose that $T\in \BS$, $T=A-B$, $A,B\in \BP$ are
such that $0\leq \tilde \phi(T)=\phi(A)-\phi(B)$. This implies
that $\phi(B)\leq \phi(A)$ which yields $B\leq A$, that is, we
have $0\leq T$. The linear transformation $\tilde \phi$ is
surjective since $\BP$ is included in its range. Moreover, it is
injective as well which follows from the fact that $\tilde \phi$
preserves the order in both directions. Now, if one further
extends $\tilde \phi$ to a linear transformation on the algebra
$B(H)$ of all bounded linear operators on $H$, one gets a linear
bijection of the $C^*$-algebra $B(H)$ which preserves the order in
both directions. Due to a well-known result of Kadison (Ref.
\cite[Corollary 5]{Kadison}) every such transformation sending the
identity to itself is a Jordan *-automorphism. Therefore, the
linear transformation
\[
A\longmapsto \sqrt{\phi(I)}^{-1}\phi(A)\sqrt{\phi(I)}^{-1}
\]
is a Jordan *-automorphism of $B(H)$. But these transformations of
$B(H)$ are well-known to be implemented by unitary-antiunitary
operators (see, for example, Ref. \cite{CVLL}). It is now easy to
infer that $\phi$ is of the desired form. This completes the proof
of the theorem.
\end{proof}

Our main result is now easy to prove.

\begin{proof}[Proof of Theorem \ref{T:lj2}]
Let $X=\phi(0)$ and consider the transformation
\[
\psi: A \mapsto \phi(A)-X.
\]
Clearly, $\psi$ is a bijection of $\BS$ preserving the order in
both directions. So, without loss of generality we can assume that
$\phi(0)=0$. Now, restricting $\phi$ onto $\BP$ we have a
bijection of $\BP$ which preserves the order in both directions.
So, we can apply Theorem~\ref{T:lj1} and obtain that there exists
an invertible bounded either linear or conjugate-linear operator
$T:H\to H$ for which we have
\begin{equation}\label{E:lj6}
\phi(A)=TAT^* \qquad (A\in \BP).
\end{equation}
It remains to show that this equality holds also for every $A\in
\BS$. Let $B\in \BS$ be arbitrary but fixed. Then there exists a
constant $K\in \R$ such that $K\leq B$ (for example, one can
choose $K=-\| B\|$). Consider the transformation
\[
A\mapsto \phi(A+K)-\phi(K)
\]
on $\BP$. Just as above, this transformation is a bijective map on $\BP$
which preserves the order in both directions. Therefore, there exists an
invertible bounded either linear or conjugate-linear operator $S:H\to H$
such that
\begin{equation}\label{E:lj2}
\phi(A+K)-\phi(K)=SAS^* \qquad (A\in \BP).
\end{equation}
If $A\geq -K, 0$, then by \eqref{E:lj6} we have
\begin{equation}\label{E:lj1}
T(A+K)T^*-\phi(K)=SAS^*.
\end{equation}
Considering this equality for another $A'$ with $A'\geq -K,0$, we see
that
\[
T(A-A')T^*=S(A-A')S^*.
\]
As the difference $A-A'$ can be an arbitrary self-adjoint operator, we
obtain that $TCT^*=SCS^*$
holds for every $C\in \BS$. It now follows from \eqref{E:lj1} that
\[
T(A+K)T^*-\phi(K)=SAS^*=TAT^*
\]
where $A\in \BS$, $A\geq -K, 0$. This yields
$\phi(K)=TKT^*=SKS^*$. We deduce from \eqref{E:lj2} that
\[
\phi(A+K)=SAS^*+\phi(K)=TAT^*+TKT^*=T(A+K)T^*
\]
holds for every $A\in \BP$. Choosing $A=B-K\geq 0$, we have
\[
\phi(B)=TBT^*.
\]
This completes the proof.
\end{proof}

We now turn to the proofs of the corollaries.

\begin{proof}[Proof of Corollary \ref{C:lj3}]
By Theorem~\ref{T:lj2} there is an invertible bounded either linear or
conjugate-linear operator $T$ on $H$ such that
\[
\phi(A)=TAT^*+\phi(0) \qquad (A\in \BS)
\]
Since $0$ is commuting with every $A\in \BS$, the same is true for
$\phi(0)$. This gives us that $\phi(0)$ is a scalar operator, that is,
there
is a $\mu \in \mathbb R$ such that $\phi(0)=\mu I$. Similarly, we have a
constant $\lambda\in \mathbb R$ such that
$TT^*=\phi(I)-\phi(0)=\lambda I$. It is trivial that $\lambda$ is
necessarily positive and then we obtain that the operator $T/\sqrt
\lambda$ is either unitary or antiunitary.
\end{proof}

\begin{proof}[Proof of Corollary \ref{C:lj4}]
It is easy to see that  $A\in \BS$ is complementary with every
$B\in \BS$ if and only if $A$ is scalar. Hence, $\phi$ preserves
the scalar operators and one can apply the argument in the proof
of Corollary~\ref{C:lj4} to get the desired form of $\phi$.
\end{proof}

In the proof of Corollary~\ref{C:lj5} we make use the following
notation. If $x,y\in H$, then $x\otimes y$ denotes the operator
defined by $(x\otimes y)z=\langle z,y\rangle x$ $(z\in H)$.

\begin{proof}[Proof of Corollary \ref{C:lj5}]
Since $0$ is the only operator in $\BS$ which is orthogonal to
every operator, we infer that $\phi(0)=0$. By Theorem~\ref{T:lj2}
we have an invertible bounded either linear or conjugate-linear
operator $T$ on $H$ such that $\phi(A)=TAT^*$ holds for every
$A\in \BS$. Without serious loss of generality we can suppose that
$T$ is linear. It now follows that for every $A,B\in \BS$ with
$AB=0$ we have $\phi(A)\phi(B)=0$ which implies that $AT^*TB=0$.
Choosing nonzero orthogonal vectors $x,y\in H$, for $A=x\otimes x$
and $B=y\otimes y$ we get $x\otimes Tx \cdot Ty \otimes y=0$ which
yields $\langle T^*Tx,y\rangle=\langle Tx,Ty\rangle=0$. So, we
have $\langle T^*Tx,y\rangle=0$ whenever $\langle x,y\rangle=0$.
This clearly implies that for every $x\in H$ there is a scalar
$\lambda_x$ such that $T^*Tx=\lambda_x x$. In another expression,
the operators $T^*T$ and $I$ are locally linearly dependent. It is
a folk result (whose proof requires only elementary linear
algebra) that in that case the operators $T^*T$ and $I$ are
necessarily linearly dependent, that is, there exists a scalar
$\lambda\in \R$ such that $T^*T=\lambda I$. Now, the proof can be
completed as in the proof of Corollary~\ref{C:lj3}.
\end{proof}

\section{Acknowledgements}
 This research was supported from the following sources:
          (1) Hungarian National Foundation for Scientific Research
          (OTKA), Grant No. T030082, T031995, (2)
          Ministry of Education, Hungary, Grant No. FKFP
          0349/2000,
          (3) Joint Hungarian-Slovene Research Project, Reg. No.
          SLO-3/00.

\newpage
\bibliographystyle{plain}

\end{document}